\documentclass{article} 
\usepackage{amsmath,amssymb,amsthm}
\usepackage{indentfirst}
\usepackage{geometry}
\usepackage{esint}
\usepackage{natbib}
\usepackage{hyperref}

\numberwithin{equation}{section}
\newtheorem{theorem}{Theorem}[section]
\newtheorem{lemma}{Lemma}[section]

\newtheorem{definition}{Definition}[section]
\newtheorem{remark}{Remark}[section]

\def\NN{\mathbb{N}}
\def\PP{\mathbb{P}}
\def\RR{\mathbb{R}}
\def\SS{\mathbb{S}}
\def\YY{\mathbb{Y}}

\title{Configuration-Dependent Lower Bounds for Approximation by Shallow ReLU$^k$ Networks on the Sphere}

\author{Tong Mao\thanks{Shenzhen Loop Area Institute, Shenzhen 518000, P. R. China}\and Jinchao Xu\footnotemark[1] \thanks{King Abdullah University of Science and Technology, Thuwal 23955, Saudi Arabia.}
}
\date{}

\begin{document}
	
\maketitle

\begin{abstract}
We establish two related but logically distinct results for shallow ReLU$^k$ neural networks on the unit sphere $\SS^d$. First, for an arbitrary set of inner neural-network parameters, the best $\mathcal{L}^2(\SS^d)$ approximation of a fixed target function with smoothness $r>\tfrac{d+2k+1}{2}$ admits an asymptotic lower bound given by a constant multiple of $n^{-1/2}\underline{h}^{k+1/2}$, where $\underline{h}$ denotes the antipodal separation distance of the normalized inner-parameter set. This lower bound depends explicitly on the parameter configuration through $\underline{h}$ and applies without additional assumptions on the parameters.

Second, for antipodally quasi-uniform parameters, $\underline{h}\simeq n^{-1/d}$, and the lower bound establishes the exact saturation order $n^{-\frac{d+2k+1}{2d}}$ for such parameter families: a target function with regularity greater than $\frac{d+2k+1}{2}$ and satisfying the required parity condition can be approximated at this rate, whereas approximation at any strictly faster rate forces the target function to be zero. Our results therefore place linearized neural-network approximation within the classical saturation framework and show that, although ReLU$^k$ network spaces can outperform finite elements of the same degree, this advantage is intrinsically limited.
\end{abstract}

\section{Introduction}

Neural networks have demonstrated remarkable approximation capabilities over the past several decades. The universal approximation theorem, established in seminal works of the late 1980s (see, e.g., \cite{cybenko1989approximation,hornik1989multilayer}), laid the theoretical foundation for their expressive power. Specifically, consider the class of shallow neural networks with a single hidden layer:
\begin{equation}
    \Sigma_n^\sigma:=\Big\{\sum_{j=1}^na_j\sigma(w_j\cdot \circ+b_j):~w_j\in\RR^d,~b_j\in\RR,~a_j\in\RR\Big\}.
\end{equation}

For suitable smooth activation functions $\sigma$, this class approximates functions in the Sobolev space $\mathcal{W}^{r,p}(\Omega)$ at the rate $\mathcal{O}(n^{-\frac{r}{d}})$ and achieves exponential convergence for analytic functions \cite{mhaskar1996neural}. Related approximation results have been established for modern architectures, including deep networks with the ReLU activation \cite{Yarotsky2017} and its higher-order variants ReLU$^k$ \cite{li2019better,he2023expressivity}, where
\begin{equation}
    \sigma_k(t)=\big(\max\{t,0\}\big)^k,\qquad t\in\RR.
\end{equation}
For instance, for any $f\in\mathcal{H}^r(\Omega)$, there exists a deep ReLU neural network $f_n$ with depth $\mathcal{O}(\log n)$ and parameter count $\mathcal{O}(n\log n)$ that achieves the approximation rate $\mathcal{O}(n^{-\frac{r}{d}})$.

For shallow networks, a central question is how the regularity of the target function influences the achievable convergence rate. Approximation by shallow ReLU$^k$ networks has been extensively studied \cite{makovoz1996random,bach2017breaking,klusowski2018approximation,xu2020finite,siegel2022sharp,mao2023rates,siegel2022optimal,mhaskar2023tractability,meng2022new}, with typical rates of the form $\mathcal{O}(n^{-\frac{d+\alpha}{2d}})$. In particular, \cite{siegel2022sharp,siegel2025optimal} established the optimal convergence rate $\mathcal{O}(n^{-\frac{d+2k+1}{2d}})$ for functions in Barron spaces $\mathcal{B}^k(\Omega)$. Building on these results, \cite{yang2024optimal,mao2024approximation} showed that, for Sobolev spaces,
\begin{equation}\label{eqn_nonlin_rate}
    \inf\limits_{f_n\in \Sigma_{n}^{\sigma_k}}\|f-f_n\|_{\mathcal{L}^2(\Omega)}\lesssim\|f\|_{\mathcal{H}^r(\Omega)}n^{-\frac{r}{d}},\quad r\leq\frac{d+2k+1}{2}.
\end{equation}

Since $\sigma_k$ is homogeneous, one may normalize $\binom{w}{b}=\theta$ so that the parameter norm is $1$, with $\theta\in\SS^d\subset\RR^{d+1}$. Recent works \cite{liu2025integral,liu2025condition} have shown that the geometry of the parameter set $\Theta=\{\theta_j^*\}_{j=1}^n\subset\SS^d$ plays an essential role in the approximation properties of shallow neural networks. Let $\rho$ be the geodesic distance on $\SS^d$, and denote the fill distance and the antipodal separation distance by
\[
h:=\max_{\eta\in\SS^d}\min_{1\leq j\leq n}\rho(\eta,\theta_j^*),\qquad \underline{h}:=\min_{i\neq j}\min\big\{\rho(\theta_i^*,\theta_j^*),\rho(\theta_i^*,-\theta_j^*)\big\}.
\]
It was proved in \cite{liu2025integral} that, after fixing the parameters $\{\theta_j^*\}_{j=1}^n$ in the nonlinear class $\Sigma_n^{\sigma_k}$ appearing in \eqref{eqn_nonlin_rate}, the linear space spanned by these $n$ fixed neurons achieves the approximation rate $\mathcal{O}(h^r)$ for functions $f\in\mathcal{H}^r(\Omega)$; see also \cite{petrushev1998approximation} for an alternative formulation. Once again, however, the regularity is restricted by the threshold $r\leq\frac{d+2k+1}{2}$.

The repeated appearance of the critical regularity
\begin{equation}
    r=\frac{d+2k+1}{2}
\end{equation}
naturally raises the questions of lower bounds and saturation. In its classical approximation-theoretic sense, saturation means more than the sharpness of an upper bound or the existence of a matching lower bound. A saturation theorem identifies a critical approximation order and contains an inverse assertion: convergence strictly faster than that order forces the target function to belong to a special exceptional class that is trivial in an appropriate sense. Thus, both the limiting rate and the characterization of the functions capable of exceeding it are essential parts of the saturation problem.

Such approximation-order barriers and saturation phenomenons are widespread in classical approximation theory. For wavelet approximation, matching lower bounds show that the convergence rate imposed by the regularity and vanishing moments of the wavelet system cannot in general be improved by assuming additional smoothness of the target function \cite{devore1992compression,cohen2001adaptive}. For the Fejér means of Fourier series, the saturation order is $n^{-\frac{1}{d}}$ when $n$ denotes the number of degrees of freedom, and convergence strictly faster than this order forces the target to be constant \cite{devore1993constructive,timan2014theory,lorentz1966approximation}. Similarly, degree-$k$ Lagrange finite element approximation has the saturation order $n^{-\frac{k+1}{d}}$: the corresponding lower bound is sharp, and convergence strictly faster than this order forces the target to belong to the exactly reproduced polynomial space $\mathbb{P}_k$ \cite{lin2014lower}.

Whether analogous lower bounds and saturation results hold for neural networks, and under what assumptions, is more subtle. The bounds in \eqref{eqn_nonlin_rate} and the earlier results \cite{siegel2022sharp,mao2023rates,mao2024approximation} suggest the exponent might be $\frac{d+2k+1}{2d}$ for shallow ReLU$^k$ networks. These upper bounds alone, however, neither provide a matching lower bound nor establish saturation in the classical sense, since they do not show that exceeding the suggested rate forces the target into a special exceptional class. Moreover, \cite{siegel2022high} showed that, in a special setting and for very smooth functions, nonlinear shallow ReLU$^k$ approximation can attain the faster rate $\mathcal{O}(n^{-(k+1)})$. The lower-bound and saturation problems for the full nonlinear class therefore remain open.

In this paper, we provide a configuration-dependent approximation lower bound for neural networks on the sphere $\SS^d$, which leads to a saturation theorem. Given parameters $\{\theta_j^*\}_{j=1}^\infty\subset\SS^d$, we consider the space
\begin{equation}\label{eqn_def_lnn_sph}
    L_{n}^k=L_{n}^k\big(\{\theta_j^*\}_{j=1}^n,\SS^d\big):=\mathrm{span}\Big\{\sigma_k(\theta_j^*\cdot\circ):~j=1,\dots,n\Big\}.
\end{equation}
For functions in $\mathcal{H}^r(\SS^d)$ satisfying the parity condition $f(\eta)=(-1)^{k+1}f(-\eta)$, the corresponding approximation upper bound is  \cite{liu2025integral}
\begin{equation}\label{eqn_lin_rate}
    \inf\limits_{f_n\in L_{n}^k}\|f-f_n\|_{\mathcal{L}^2(\SS^d)}\lesssim h^r\|f\|_{\mathcal{H}^r(\SS^d)},\quad r\leq\frac{d+2k+1}{2}.
\end{equation}
When $\{\theta_j^*\}_{j=1}^n$ is quasi-uniform, \eqref{eqn_lin_rate} yields the sharp rate $\mathcal{O}(n^{-\frac rd})$.

Our first result is a configuration-dependent lower bound. For an arbitrary parameter set $\{\theta_j^*\}_{j=1}^n\subset\SS^d$ and any $f\in\mathcal{H}^r(\SS^d)$ with $r>\frac{d+2k+1}{2}$, we prove
\begin{equation}
    \inf\limits_{f_n\in L_n^k}\|f-f_n\|_{\mathcal{L}^2(\SS^d)}\gtrsim n^{-\frac12}\underline{h}^{k+\frac12}\|f\|_{\mathcal{L}^2(\SS^d)}.
\end{equation}
This estimate does not require quasi-uniformity or $\underline{h}>0$ and records explicitly, through $\underline{h}$, the effect of near-coincident and near-antipodal parameters.

The saturation conclusion is obtained after imposing antipodal quasi-uniformity on the parameter families. If $h\lesssim\underline{h}$, standard covering and packing estimates give $\underline{h}\simeq n^{-1/d}$, and the configuration-dependent lower bound becomes
\begin{equation}\label{eqn_satu_lower}
    \inf\limits_{f_n\in L_n^k}\|f-f_n\|_{\mathcal{L}^2(\SS^d)}\gtrsim n^{-\frac{d+2k+1}{2d}}\|f\|_{\mathcal{L}^2(\SS^d)}.
\end{equation}
On the class
\begin{equation}
    \mathcal{F}^{d,k}_*:=\left\{f\in\mathcal{H}^r(\SS^d):~r>\frac{d+2k+1}{2},~f(\eta)=(-1)^{k+1}f(-\eta)\right\},
\end{equation}
the upper estimate \eqref{eqn_lin_rate} and the lower estimate \eqref{eqn_satu_lower} together yield the saturation phenomenon:
\begin{enumerate}
    \item Any $f\in\mathcal{F}^{d,k}_*$ can be approximated with rate $\mathcal{O}(n^{-\frac{d+2k+1}{2d}})$;
    \item $f\in\mathcal{F}^{d,k}_*$ can be approximated with rate $o(n^{-\frac{d+2k+1}{2d}})$ if and only if $f\equiv0$.
\end{enumerate}
The first statement identifies the attainable critical rate, while the second is the inverse part of the saturation theorem: within $\mathcal{F}^{d,k}_*$, exceeding that rate forces the target to be trivial. Antipodal quasi-uniformity therefore converts the general $\underline{h}$-dependent lower bound into a sharp rate in terms of $n$ and, together with the matching upper bound, establishes saturation. It is not an assumption of the general lower bound itself.

The proof reflects the antipodal geometry of ReLU$^k$ ridge functions. We decompose the squared norm into dyadic blocks indexed by the actual spherical-harmonic degree $2m+I_k$. A parity projection of the associated filtered Jacobi kernels produces localization at both $t=1$ and $t=-1$, making $\underline{h}$ the natural geometric scale in the matrix lower bound.

To the best of our knowledge, this is the first saturation result of this type for linearized shallow ReLU$^k$ networks. It complements the upper bounds of \cite{liu2025integral} beyond the critical regularity $\frac{d+2k+1}{2}$ and quantifies the limiting approximation behavior of arbitrary parameter configurations through $\underline{h}$. In the antipodally quasi-uniform case, the saturation rate $n^{-\frac{d+2k+1}{2d}}$ also shows that the improvement of these linearized neural-network spaces over degree-$k$ finite elements, whose saturation order is $\frac{k+1}{d}$, cannot continue beyond the critical Sobolev smoothness.

Throughout this paper, following \cite{xu1992iterative}, the notation $f(x)\lesssim B$ means that $f(x)\leq Cg(x)$ for a constant $C>0$ independent of the discretization parameters under consideration. We write $f(x)\gtrsim g(x)$ if $g(x)\lesssim f(x)$ and $f(x)\simeq g(x)$ if both $f(x)\lesssim g(x)$ and $f(x)\gtrsim g(x)$.

\section{Localized spherical polynomials}

In this section, we establish a decomposition of the squared $\mathcal{L}^2$ norm for functions in $L_n^k$. Specifically, we construct matrices
\[
Q_q=\big(L_q(\theta_i^*\cdot\theta_j^*)\big)_{i,j=1}^n,\qquad q=0,1,\dots,
\]
such that every $f_n=\sum_{j=1}^na_j\sigma_k(\theta_j^*\cdot\circ)\in L_n^k$ satisfies
\[
\|f_n\|_{\mathcal{L}^2(\SS^d)}^2=\sum_{m=0}^k\widehat{\sigma_k}(m)^2a^\top P(m)a+\sum_{q=0}^\infty a^\top Q_q a,
\]
where $a=(a_1,\ldots,a_n)^\top$. The first term is the finite low-degree polynomial contribution; the dyadic matrices describe the remaining harmonic degrees.

This decomposition plays a central role in establishing the saturation phenomenon. For smooth functions $f\in\mathcal{H}^r(\SS^d)$, classical polynomial approximation theory controls the decay of their high-degree components. For functions in $L_n^k$, the spectral properties of the matrices $Q_q$ instead yield a configuration-dependent lower bound involving $n^{-1/2}\underline{h}^{k+1/2}$; the sharp power $n^{-\frac{d+2k+1}{2d}}$ follows in the antipodally quasi-uniform case.

The key to establishing this lower bound lies in showing that the matrices $Q_q$ are strongly diagonally dominant. This property emerges from the localization characteristics of spherical harmonic polynomials—a fundamental concept in approximation theory that has been extensively studied \cite{petrushev2005localized,ivanov2010sub,dai2013approximation,xu2024highly}. The localization ensures that the influence of each basis function remains concentrated, leading to the diagonal dominance that ultimately constrains the approximation power of linearized neural networks.

\subsection{Spherical harmonics and Legendre polynomials}
We begin with some standard notation together with basic facts from spherical harmonic analysis. The material in this subsection is classical and can be found in \cite{liu2025integral,dai2013approximation}. We write $\NN=\{1,2,\ldots\}$, $\NN_0=\{0,1,2,\ldots\}$, and $\SS^d:=\{\eta\in\RR^{d+1}:|\eta|=1\}$. For $\eta,\theta\in\SS^d$, the geodesic distance is
\[
\rho(\eta,\theta):=\arccos(\eta\cdot\theta).
\]
Let $\omega_d:=\int_{\SS^d}1\,d\eta$ denote the surface area of $\SS^d$. We use the normalized surface measure
\[
\fint_{\SS^d} f(\eta)\,d\eta:=\frac{1}{\omega_d}\int_{\SS^d} f(\eta)\,d\eta,\qquad f\in\mathcal{L}^1(\SS^d),
\]
and the induced $\mathcal{L}^2$ inner product and norm
\[
\langle f,g\rangle_{\mathcal{L}^2(\SS^d)}:=\fint_{\SS^d}f(\eta)g(\eta)\,d\eta,\qquad \|f\|_{\mathcal{L}^2(\SS^d)}^2:=\langle f,f\rangle_{\mathcal{L}^2(\SS^d)}.
\]

Let $\PP_m(\SS^d)$ be the space obtained by restricting to $\SS^d$ all polynomials in $\RR^{d+1}$ of total degree at most $m$. Its dimension is
\[
\dim\PP_m(\SS^d)=
\begin{cases}
\binom{d+1+m}{m}, & m=0,1,\\[6pt]
\binom{d+1+m}{m}-\binom{d-1+m}{m-2}, & m\ge 2.
\end{cases}
\]
Let $\YY_m$ be the orthogonal complement of $\PP_{m-1}(\SS^d)$ in $\PP_m(\SS^d)$. This is the space of spherical harmonics of degree $m$. Fix a real $\mathcal{L}^2$-orthonormal basis $\{Y_{m,\ell}\}_{\ell=1}^{N(m)}\subset\YY_m$. Its dimension is
\[
N(0)=1,\qquad
N(m)=\frac{2m+d-1}{m}\binom{m+d-2}{d-1},\quad m\ge1.
\]

Every $f\in\mathcal{L}^2(\SS^d)$ admits the harmonic expansion
\[
f(\eta)=\sum_{m=0}^\infty\sum_{\ell=1}^{N(m)}\widehat f(m,\ell)Y_{m,\ell}(\eta),\qquad \widehat f(m,\ell):=\langle f,Y_{m,\ell}\rangle_{\mathcal{L}^2(\SS^d)},
\]
and the $\mathcal{L}^2$ projection onto $\YY_m$ is
\[
\Pi_m f:=\sum_{\ell=1}^{N(m)}\widehat f(m,\ell)Y_{m,\ell}.
\]
We also write
\[
\mathcal{P}_M:=\sum_{m=0}^M\Pi_m
\]
for the $\mathcal{L}^2$-orthogonal projection onto $\PP_M(\SS^d)$.
Then Parseval’s identity
\[
\|f\|_{\mathcal{L}^2(\SS^d)}^2=\sum_{m=0}^\infty\sum_{\ell=1}^{N(m)}|\widehat f(m,\ell)|^2=\sum_{m=0}^{\infty}\|\Pi_m f\|_{\mathcal{L}^2(\SS^d)}^2
\]
gives the standard definition of Sobolev spaces.
\begin{definition}[Sobolev spaces on the sphere]\label{def:sobolev_sphere}
    For $r>0$, the Sobolev space $\mathcal{H}^r(\SS^d)$ is defined as $\mathcal{H}^r(\SS^d)=\{f\in \mathcal{L}^2(\SS^d):~\|f\|_{\mathcal{H}^r(\SS^d)}<\infty\}$, with norm squared
    \begin{equation}\label{eqn:Sob_norm_Parseval}
        \|f\|_{\mathcal{H}^r(\SS^d)}^2=\|f\|_{\mathcal{L}^2(\SS^d)}^2+\sum_{m=1}^\infty m^{2r}\|\Pi_m f\|_{\mathcal{L}^2(\SS^d)}^2 = \sum_{m=0}^\infty\sum_{\ell=1}^{N(m)}(m^{2r}+1)|\widehat f(m,\ell)|^2.
    \end{equation}
\end{definition}

Define the space $\mathcal{L}^2_{w_d}([-1,1])$ by
\begin{equation}
    \langle f,g\rangle_{w_d}:=\int_{-1}^1f(t)g(t)(1-t^2)^{\frac{d-2}{2}}\,dt,\qquad\|f\|_{\mathcal{L}^2_{w_d}([-1,1])}:=\langle f,f\rangle_{w_d}^{\frac{1}{2}}.
\end{equation}
We use the following Legendre polynomials as an orthogonal basis of this space (see, e.g., \cite{szego1975orthogonal}):
\begin{equation}\label{eqn:lambdam_pm}
    p_m(t)=\lambda_m(1-t^2)^{-\frac{d-2}{2}}\left(\frac{d}{dt}\right)^m\left[(1-t^2)^{m+\frac{d-2}{2}}\right],\qquad t\in[-1,1],
\end{equation}
where
\begin{equation*}
    \lambda_m=\frac{\omega_{d}}{\omega_{d-1}}\frac{N(m)}{\Gamma(m+d/2)}\sqrt{\frac{(2m+d-1)\Gamma(m+d-1)}{2^{2m+d-1}\Gamma(m+1)}},\qquad m\in\NN_0,
\end{equation*}
are chosen such that
\begin{equation}\label{eqn:sum_Y_nl}
p_m(\eta\cdot \theta) = \sum_{\ell=1}^{N(m)} Y_{m,\ell}(\eta)Y_{m,\ell}(\theta).
\end{equation}

The function $\sigma_k\in \mathcal{L}^2_{w_d}([-1,1])$ has the Legendre expansion
\begin{equation}\label{eqn:Legendre_expansion_ReLUk}
\sigma_k=\sum_{m=0}^\infty\widehat{\sigma_k}(m)p_m,
\end{equation}
where the Legendre coefficients are given by
\[
\widehat{\sigma_k}(m)=\frac{\langle p_m,\sigma_k\rangle_{w_d}}{\|p_m\|_{\mathcal{L}^2_{w_d}([-1,1])}^2}.
\]
Define
\begin{equation}
    E_{\sigma_k}:=\left\{m\in\NN_0:~\widehat{\sigma_k}(m)\neq0\right\}.
\end{equation}
Then, by \cite[Appendix D.2]{bach2017breaking},
\begin{equation}\label{eqn:hat_sigma_large}
    \begin{split}
        &E_{\sigma_k}=\left\{m\geq k+1:~m-k\text{ is odd}\right\}\cup\{0,\dots,k\},\\
        &\widehat{\sigma_k}(m)=\frac{\omega_{d-1}k!\Gamma(d/2)}{\omega_d}\frac{(-1)^{(m-k-1)/2}\Gamma(m-k)}{2^m\Gamma\left(\frac{m-k+1}{2}\right)\Gamma\left(\frac{m+d+k+1}{2}\right)},\qquad m\geq k+1,\quad m-k\text{ is odd}.
    \end{split}
\end{equation}

Finally, for notational simplicity, we follow \cite[Lemma 3]{liu2025integral} and set
\begin{equation}\label{eqn:theta}
    \xi(t)=\left(\frac{\omega_{d-1}}{\omega_d}\frac{k!\Gamma(d/2)}{2^{k+1}\sqrt{\pi}}\right)^2 \left(\frac{\Gamma\left(\frac{t-k}{2}\right)}{\Gamma\left(\frac{t+d+k+1}{2}\right)}\right)^2.
\end{equation}
By the duplication formula,
\begin{equation}\label{eqn:hat_sigma_xi}
    \widehat{\sigma_k}(m)^2=\xi(m),\qquad m\geq k+1,\quad m-k\text{ is odd}.
\end{equation}

\subsection{Representing the norm of \texorpdfstring{$f_n$}{f\_n}}
In this subsection, we derive an explicit representation of the norm of $f_n$ in terms of Legendre polynomials and spherical harmonics. We first introduce the notation
\begin{equation}
    I_k=\begin{cases}
        0,&k\text{ is odd},\\
        1,&k\text{ is even}.
    \end{cases}
\end{equation}

As in \cite{liu2025integral}, the norm of the function
\[
f_n(\eta)=\sum_{j=1}^na_j\sigma_k(\theta_j^*\cdot\eta)
\]
can be written as
\begin{equation}\label{eqn_fn_norm_explicit}
    \begin{split}
        \|f_n\|_{\mathcal{L}^2(\SS^d)}^2=&\Big\|\sum_{j=1}^na_j\sum_{m=0}^{\infty}\widehat{\sigma_k}(m)\sum_{\ell=1}^{N(m)}Y_{m,\ell}(\theta_j^*)Y_{m,\ell}\Big\|_{\mathcal{L}^2(\SS^d)}^2=\sum_{m=0}^{\infty}\widehat{\sigma_k}(m)^2\sum_{\ell=1}^{N(m)}\Big(\sum_{j=1}^na_jY_{m,\ell}(\theta_j^*)\Big)^2\\
        =&\sum_{m=0}^{\infty}\widehat{\sigma_k}(m)^2a^\top P(m)a,
    \end{split}
\end{equation}
where
\[
P(m)=\Big(\sum_{\ell=1}^{N(m)}Y_{m,\ell}(\theta_i^*)Y_{m,\ell}(\theta_j^*)\Big)_{i,j=1}^n=\big(p_m(\theta_i^*\cdot\theta_j^*)\big)_{i,j=1}^n.
\]

It is known (see, e.g., \cite[(3.6)]{petrushev2005localized}) that there exists a smooth function $\zeta$ satisfying
\begin{equation*}
\begin{aligned}
&\zeta\in\mathcal{C}^\infty(\RR),\qquad &&\zeta\geq0,\qquad \operatorname{supp}(\zeta)\subset[1/2,2],\\
&\zeta(t)>c_1>0,\qquad &&t\in[3/5,5/3],\\
&\zeta(t)+\zeta(2t)=1,\qquad &&t\in[1/2,1].
\end{aligned}
\end{equation*}
In particular, $0\leq\zeta\leq1$ and
\[
1=\sum_{q=0}^\infty\zeta(2^{-q}\ell),\qquad \ell\in\NN.
\]
Then we can write \eqref{eqn_fn_norm_explicit} as
\begin{equation}\label{eqn_fn_Qq_explicit}
    \begin{split}
        \|f_n\|_{\mathcal{L}^2(\SS^d)}^2
        =&\sum_{m=0}^k\widehat{\sigma_k}(m)^2a^\top P(m)a+\sum_{\substack{m\geq0\\2m+I_k\geq k+1}}\xi(2m+I_k)a^\top P(2m+I_k)a\\
        =&\sum_{m=0}^k\widehat{\sigma_k}(m)^2a^\top P(m)a+\sum_{q=0}^\infty a^\top Q_q a,
    \end{split}
\end{equation}
where
\begin{equation}\label{eqn_def_Qq}
Q_q=\sum_{\substack{m\geq0\\2m+I_k\geq k+1}}\zeta\left(\frac{2m+I_k}{2^q}\right)\xi(2m+I_k)P(2m+I_k).
\end{equation}

Moreover, let $q\in\NN_0$ and suppose that $2^q\geq k+1$. Since $\mathcal{P}_{2^q}$ is the orthogonal projection onto $\PP_{2^q}(\SS^d)$ and the harmonic-degree support of $Q_{q+1}$ is contained in $[2^q+1,2^{q+2}-1]$, we have
\begin{equation}\label{eqn_fn-Pfn}
    \begin{split}
        \big\|\mathcal{P}_{2^{q+2}}(f_n)-\mathcal{P}_{2^{q}}(f_n)\big\|_{\mathcal{L}^2(\SS^d)}^2
        =&\sum_{m=2^q+1}^{2^{q+2}}\widehat{\sigma_k}(m)^2a^\top P(m)a\\
        \geq&\sum_{\substack{m\geq0\\2m+I_k\geq k+1}}\zeta\left(\frac{2m+I_k}{2^{q+1}}\right)\xi(2m+I_k)a^\top P(2m+I_k)a
        =a^\top Q_{q+1}a.
    \end{split}
\end{equation}

The estimates below are stated for arbitrary parameter sets. Antipodal quasi-uniformity is used only when the final $\underline{h}$-dependent estimate is converted into an explicit rate in terms of $n$.

\subsection{Summation of Jacobi polynomials and localization}

The use of summability kernels to obtain localization goes back at least to the work of Fejér \cite{fejer1903untersuchungen}. Smoothly filtered analogues were subsequently developed for trigonometric, Jacobi, and spherical polynomial expansions; see, for example, \cite{mhaskarPrestin2005local,mhaskar2004polynomialLocal,mhaskar2005sphereLocalized}. For Jacobi expansions, the localized-kernel and frame construction in \cite{petrushev2005localized} provides the estimate needed below; related exponentially and sub-exponentially localized constructions can be found in \cite{filbir2009filter,ivanov2010sub}.

The full filtered Jacobi kernel is localized near $t=1$. The sums arising in our analysis, however, contain only even or only odd degrees and are therefore localized near both $t=1$ and $t=-1$. We obtain this two-endpoint localization by applying the parity projector to the full kernel. Throughout this subsection, we assume $d\geq2$ and write $\|g\|_{w_d}:=\|g\|_{\mathcal{L}^2_{w_d}([-1,1])}$.

\begin{theorem}[Theorem~2.4, \cite{petrushev2005localized}]\label{thm_petrushev_xu_localized}
Fix $\gamma\in\NN$, and let $A\in\mathcal{C}^{\gamma+2d-1}(\RR)$ be supported in $[1/2,2]$.
For $q\in\NN_0$, let
\begin{equation*}
    K_q(t):=\sum_{m=0}^\infty 
    A\left(\frac{m}{2^q}\right)
    \frac{p_m(1)p_m(t)}{\|p_m\|_{w_d}^2},
    \qquad -1\le t\le1.
\end{equation*}
Then
\begin{equation*}
    |K_q(t)|
    \lesssim
    \frac{2^{q\frac{d+1}{2}}\max_{0\le\nu\le\gamma+2d-1}
    \|A^{(\nu)}\|_{\mathcal{L}^1([1/2,2])}}
    {(\sqrt{1+t}+2^{-q})^{\frac{d-1}{2}}
     (\sqrt{1-t}+2^{-q})^{\frac{d-1}{2}}
     (1+2^q\sqrt{1-t})^\gamma}.
\end{equation*}
The implied constant depends only on $d$ and $\gamma$.
\end{theorem}

The particular support interval $[1/2,2]$ only fixes the spectral scale. The same estimate holds when $A$ is supported in any fixed interval $[a,b]\subset(0,\infty)$, with the implied constant additionally depending on $a$ and $b$. Indeed, one may decompose $A$ into finitely many smooth pieces supported in dyadic intervals of the form $2^j[1/2,2]$, rescale each piece, and apply Theorem~\ref{thm_petrushev_xu_localized} at the corresponding dyadic level.

\begin{theorem}\label{thm:localize2}
Let $K\in\NN$ and $q\in\NN_0$, and let $\varphi\in\mathcal{C}^{K+2d-1}([0,\infty))$ satisfy $\operatorname{supp}(\varphi)\subset[1/2,2]$. Define
\begin{equation}\label{eqn:def_L_localized}
L(t)=\sum_{m=0}^\infty\varphi\left(\frac{2m+I_k}{2^q}\right)p_{2m+I_k}(t).
\end{equation}
Then, for every $t\in[-1,1]$,
\begin{equation}\label{eqn:L_localized}
|L(t)|\lesssim\max_{0\leq\beta\leq K+2d-1}\|\varphi^{(\beta)}\|_{\mathcal{L}^1([0,\infty))}\frac{2^{qd}}{(1+2^q\sqrt{1-t^2})^K}.
\end{equation}
The implied constant depends only on $d$ and $K$.
\end{theorem}

\begin{proof}
Extend $\varphi$ by zero to $\RR$ and define
\begin{equation*}
\mathcal K_q(t):=\sum_{M=0}^\infty\varphi\left(\frac{M}{2^q}\right)\frac{p_M(1)p_M(t)}{\|p_M\|_{w_d}^2}.
\end{equation*}
Under the normalization used above,
\begin{equation*}
p_M(1)=N(M),\qquad \|p_M\|_{w_d}^2=\frac{\omega_d}{\omega_{d-1}}N(M),
\end{equation*}
and hence
\begin{equation*}
\frac{p_M(1)p_M(t)}{\|p_M\|_{w_d}^2}=\frac{\omega_{d-1}}{\omega_d}p_M(t).
\end{equation*}
Applying Theorem~\ref{thm_petrushev_xu_localized} with $\gamma=K$ and observing that
\begin{equation*}
(\sqrt{1+t}+2^{-q})(\sqrt{1-t}+2^{-q})\geq2^{-q}\big(\sqrt{1+t}+\sqrt{1-t}\big)\geq\sqrt{2}\,2^{-q},
\end{equation*}
we obtain
\begin{equation}\label{eqn:full_kernel_one_endpoint}
|\mathcal K_q(t)|\lesssim\max_{0\leq\beta\leq K+2d-1}\|\varphi^{(\beta)}\|_{\mathcal{L}^1([0,\infty))}\frac{2^{qd}}{(1+2^q\sqrt{1-t})^K}.
\end{equation}

Since $p_M(-t)=(-1)^M p_M(t)$, parity projection gives
\begin{equation*}
\begin{split}
\frac{\mathcal K_q(t)+(-1)^{I_k}\mathcal K_q(-t)}{2}
&=\frac{\omega_{d-1}}{\omega_d}\sum_{\substack{M\geq0\\ M\equiv I_k\pmod 2}}\varphi\left(\frac{M}{2^q}\right)p_M(t)\\
&=\frac{\omega_{d-1}}{\omega_d}\sum_{m=0}^\infty\varphi\left(\frac{2m+I_k}{2^q}\right)p_{2m+I_k}(t)
=\frac{\omega_{d-1}}{\omega_d}L(t).
\end{split}
\end{equation*}
Therefore, applying \eqref{eqn:full_kernel_one_endpoint} at $t$ and $-t$ yields
\begin{equation*}
|L(t)|\lesssim\max_{0\leq\beta\leq K+2d-1}\|\varphi^{(\beta)}\|_{\mathcal{L}^1([0,\infty))}2^{qd}\left(\frac{1}{(1+2^q\sqrt{1-t})^K}+\frac{1}{(1+2^q\sqrt{1+t})^K}\right).
\end{equation*}
Finally, an elementary estimation 
\begin{equation*}
\frac{1}{(1+2^q\sqrt{1-t})^K}+\frac{1}{(1+2^q\sqrt{1+t})^K}\lesssim\frac{1}{(1+2^q\sqrt{1-t^2})^K}
\end{equation*}
proves \eqref{eqn:L_localized}.
\end{proof}

\subsection{Lower bounds for the matrices \texorpdfstring{$Q_q$}{Q\_q}}

In this subsection, we establish a lower bound for $Q_q$ by combining the two-endpoint localization of $L_q$ with the antipodal separation of the parameters.

\begin{lemma}\label{lem_Q_q_norm}
Let $K\in\NN$ with $K>d$, let $\{Q_q\}_{q=0}^\infty$ be defined by \eqref{eqn_def_Qq}, and let
\begin{equation*}
\underline{h}:=\min_{i\neq j}\min\{\rho(\theta_i^*,\theta_j^*),\rho(\theta_i^*,-\theta_j^*)\}>0.
\end{equation*}
Then, for all $q\geq1+\log_2(k+1)$,
\begin{equation}\label{eqn:Qq_off_diagonal}
\max_{1\leq j\leq n}\sum_{i\neq j}|(Q_q)_{i,j}|\lesssim2^{-q(2k+1+K)}\underline{h}^{-K}.
\end{equation}
Moreover, there exists a constant $C_3>0$ such that
\begin{equation}\label{eqn:Qq_lower_bound}
Q_q\gtrsim2^{-q(2k+1)}I_{n\times n},\qquad q\geq\log_2\left(\frac{C_3}{\underline{h}}\right).
\end{equation}
Consequently, for every $f_n\in L_n^k$, represented by the coefficient vector $a=(a_1,\ldots,a_n)^\top$,
    \begin{equation}
        \big\|\mathcal{P}_{2^{q+2}}(f_n)-\mathcal{P}_{2^{q}}(f_n)\big\|_{\mathcal{L}^2(\SS^d)}\gtrsim2^{-q(k+\frac12)}\|a\|_2,\qquad q\geq\log_2\left(\frac{C_3}{\underline{h}}\right).
    \end{equation}
The constants are independent of $q$, $n$, and the $\{\theta_j^*\}_{j=1}^n$.
\end{lemma}

\begin{proof}
Since $2^{q-1}>k$, on the support of $\zeta((2m+I_k)/2^q)$, $2m+I_k\geq k+1$. So the lower-frequency restriction in \eqref{eqn_def_Qq} is automatic.
By definition,
\begin{equation*}
Q_q=\big(L_q(\theta_i^*\cdot\theta_j^*)\big)_{i,j=1}^n,
\end{equation*}
where
\begin{equation}\label{eqn:def_Lq}
L_q(t)=\sum_{m=0}^\infty\varphi_q\left(\frac{2m+I_k}{2^q}\right)p_{2m+I_k}(t),\qquad \varphi_q(t):=\zeta(t)\xi(2^qt).
\end{equation}
Indeed,
\begin{equation*}
\varphi_q\left(\frac{2m+I_k}{2^q}\right)=\zeta\left(\frac{2m+I_k}{2^q}\right)\xi(2m+I_k).
\end{equation*}

By \cite[(3.21)]{liu2025integral}, for every fixed $\nu\geq0$,
\begin{equation*}
|\xi^{(\nu)}(s)|\lesssim s^{-d-2k-1-\nu}
\end{equation*}
for all sufficiently large $s$. Hence, for sufficiently large $q$ and $0\leq\beta\leq K+2d-1$, Leibniz's rule gives
\begin{equation}\label{eqn:varphi_q_derivative}
\begin{split}
\|\varphi_q^{(\beta)}\|_{\mathcal{L}^1([0,\infty))}
&\leq\sum_{\nu=0}^\beta\binom{\beta}{\nu}2^{q\nu}\int_{1/2}^2|\zeta^{(\beta-\nu)}(t)|\,|\xi^{(\nu)}(2^qt)|\,dt\\
&\lesssim\sum_{\nu=0}^\beta2^{q\nu}2^{-q(d+2k+1+\nu)}\int_{1/2}^2|\zeta^{(\beta-\nu)}(t)|t^{-d-2k-1-\nu}\,dt\\
&\lesssim2^{-q(d+2k+1)}.
\end{split}
\end{equation}
Theorem~\ref{thm:localize2} therefore implies
\begin{equation}\label{eqn_Lq_esti}
|L_q(t)|\lesssim\frac{2^{-q(2k+1)}}{(1+2^q\sqrt{1-t^2})^K},\qquad -1\leq t\leq1.
\end{equation}

For $i\neq j$, we have
\begin{equation*}
\sqrt{1-(\theta_i^*\cdot\theta_j^*)^2}
=\sin\left(\min\{\rho(\theta_i^*,\theta_j^*),\rho(\theta_i^*,-\theta_j^*)\}\right)
\gtrsim\min\{\rho(\theta_i^*,\theta_j^*),\rho(\theta_i^*,-\theta_j^*)\}.
\end{equation*}
For each fixed $j$, the separation condition and the standard spherical packing estimate give
\begin{equation*}
\#\left\{i\neq j:2^p\underline{h}\leq\min\{\rho(\theta_i^*,\theta_j^*),\rho(\theta_i^*,-\theta_j^*)\}<2^{p+1}\underline{h}\right\}\lesssim2^{pd}.
\end{equation*}
Consequently, by \eqref{eqn_Lq_esti},
\begin{equation*}
\begin{split}
\sum_{i\neq j}|(Q_q)_{i,j}|
&\lesssim2^{-q(2k+1)}\sum_{p=0}^\infty\frac{2^{pd}}{(1+2^{q+p}\underline{h})^K}\\
&\lesssim2^{-q(2k+1)}(2^q\underline{h})^{-K}\sum_{p=0}^\infty2^{-p(K-d)}
\lesssim2^{-q(2k+1+K)}\underline{h}^{-K},
\end{split}
\end{equation*}
where the series converges because $K>d$. This proves \eqref{eqn:Qq_off_diagonal}.

For the diagonal entries, using $p_M(1)=N(M)$, we have
\begin{equation*}
L_q(1)=\sum_{m=0}^\infty\zeta\left(\frac{2m+I_k}{2^q}\right)\xi(2m+I_k)N(2m+I_k).
\end{equation*}
The Gamma-ratio asymptotics in \eqref{eqn:theta} yield
\begin{equation*}
\xi(M)\simeq M^{-d-2k-1},\qquad N(M)\simeq M^{d-1},
\end{equation*}
for sufficiently large $M$. Since $\zeta\geq c_1>0$ on $[3/5,5/3]$ and all the summands are nonnegative,
\begin{equation}\label{eqn:Lq_diagonal}
\begin{split}
L_q(1)
&\geq c\sum_{\substack{m\geq0\\ (3/5)2^q\leq2m+I_k\leq(5/3)2^q}}\xi(2m+I_k)N(2m+I_k)\\
&\gtrsim\sum_{\substack{m\geq0\\ (3/5)2^q\leq2m+I_k\leq(5/3)2^q}}(2m+I_k)^{-2k-2}
\gtrsim2^{-q(2k+1)}.
\end{split}
\end{equation}

Combining \eqref{eqn:Qq_off_diagonal} and \eqref{eqn:Lq_diagonal}, there exists $C_3>0$ such that
\begin{equation*}
\max_{1\leq j\leq n}\sum_{i\neq j}|(Q_q)_{i,j}|\leq\frac12L_q(1)
\end{equation*}
whenever $2^q\underline{h}\geq C_3$. Enlarging $C_3$ if necessary also ensures that $q$ is sufficiently large for the preceding estimates. Since $Q_q$ is symmetric, for every $a\in\RR^n$,
\begin{equation*}
a^\top Q_q a\geq\sum_{i=1}^n\left(L_q(1)-\sum_{j\neq i}|(Q_q)_{i,j}|\right)a_i^2\geq\frac12L_q(1)\|a\|_2^2\gtrsim2^{-q(2k+1)}\|a\|_2^2.
\end{equation*}
This proves \eqref{eqn:Qq_lower_bound}.
Finally, by \eqref{eqn_fn-Pfn} and Lemma~\ref{lem_Q_q_norm},
\begin{equation*}
\begin{split}
    \big\|f_n-\mathcal{P}_{2^{q}}(f_n)\big\|_{\mathcal{L}^2(\SS^d)}^2
    \geq a^\top Q_{q+1}a
    \gtrsim2^{-(q+1)(2k+1)}\|a\|_2^2\simeq2^{-q(2k+1)}\|a\|_2^2.
\end{split}
\end{equation*}
\end{proof}

\section{Saturation phenomenon for linearized \texorpdfstring{ReLU$^k$}{ReLU-k} neural networks}

We now establish a configuration-dependent lower bound for linearized ReLU$^k$ approximation. The general estimate is expressed in terms of $\underline{h}$ and does not require quasi-uniformity. The saturation order $\frac{d+2k+1}{2d}$ follows as the sharp $n$-dependent consequence for antipodally quasi-uniform center families.

\begin{theorem}\label{thm_main}
    Let $d,k\in\NN$, let $\{\theta_j^*\}_{j=1}^n\subset\SS^d$ satisfy $\underline{h}>0$, and suppose that $r>\frac{d+2k+1}{2}$. Then, for every $f\in\mathcal{H}^r(\SS^d)$,
    \begin{equation}
        \inf\limits_{f_n\in L_n^k}\|f-f_n\|_{\mathcal{L}^2(\SS^d)}\gtrsim n^{-\frac12}\underline{h}^{k+\frac12}\|f\|_{\mathcal{L}^2(\SS^d)}.
    \end{equation}
    In particular, if $\{\theta_j^*\}_{j=1}^n$ is antipodally quasi-uniform,
    \begin{equation}
        \inf\limits_{f_n\in L_n^k}\|f-f_n\|_{\mathcal{L}^2(\SS^d)}\gtrsim n^{-\frac{d+2k+1}{2d}}\|f\|_{\mathcal{L}^2(\SS^d)},
    \end{equation}
    where the corresponding constant is independent of $n$.
\end{theorem}
\begin{proof}
    Without loss of generality, assume $\|f_n-f\|_{\mathcal{L}^2(\SS^d)}\leq(1-\frac{1}{\sqrt{2}})\|f\|_{\mathcal{L}^2(\SS^d)}$. Then we have $\|f_n\|_{\mathcal{L}^2(\SS^d)}\geq\frac{1}{\sqrt{2}}\|f\|_{\mathcal{L}^2(\SS^d)}$ and
    \begin{equation*}
        \begin{split}
            n\|a\|_2^2\geq&\sup\limits_{\eta\in\SS^d}\Big(n\sum_{j=1}^na_j^2\Big)\Big(\frac{1}{n}\sum_{j=1}^n\sigma_k(\theta_j^*\cdot\eta)^2\Big)\geq\sup\limits_{\eta\in\SS^d}\Big(\sum_{j=1}^na_j\sigma_k(\theta_j^*\cdot\eta)\Big)^2\\
            \geq&\|f_n\|_{\mathcal{L}^2(\SS^d)}^2\geq\frac{\|f\|_{\mathcal{L}^2(\SS^d)}^2}{2},
        \end{split}
    \end{equation*}
    which implies
    \begin{equation}\label{eqn:est_a_lower}
    \|a\|_2^2\gtrsim n^{-1}\|f\|_{\mathcal{L}^2(\SS^d)}^2.
\end{equation}

Let $q=\lceil\log_2\left(\frac{C_3}{\underline{h}}\right)\rceil$, classical approximation theory gives (see, e.g., \cite{dai2013approximation,devore1993constructive})
\begin{equation*}
    \|f-\mathcal{P}_{2^{q}}(f)\|_{\mathcal{L}^2(\SS^d)}\lesssim2^{-q r}\|f\|_{\mathcal{H}^r(\SS^d)}\simeq \underline{h}^r\|f\|_{\mathcal{H}^r(\SS^d)}.
\end{equation*}
By Lemma~\ref{lem_Q_q_norm},
\begin{equation*}
\begin{split}
    \big\|f_n-\mathcal{P}_{2^{q}}(f_n)\big\|_{\mathcal{L}^2(\SS^d)}^2\gtrsim2^{-q(2k+1)}\|a\|_2^2
    \simeq n^{-1}\underline{h}^{2k+1}\|f\|_{\mathcal{L}^2(\SS^d)}^2.
\end{split}
\end{equation*}
Writing
\begin{equation}\label{eqn_est_fn-f}
    \begin{split}
        &\big\|f_n-f\big\|_{\mathcal{L}^2(\SS^d)}=\Big(\big\|f_n-\mathcal{P}_{2^{q}}(f_n)-(f-\mathcal{P}_{2^{q}}(f))\big\|_{\mathcal{L}^2(\SS^d)}^2+\big\|\mathcal{P}_{2^{q}}(f_n-f)\big\|_{\mathcal{L}^2(\SS^d)}^2\Big)^{\frac{1}{2}}\\
        \geq&\big\|f_n-\mathcal{P}_{2^{q}}(f_n)-(f-\mathcal{P}_{2^{q}}(f))\big\|_{\mathcal{L}^2(\SS^d)},
    \end{split}
\end{equation}
then
\begin{equation}
    \|f_n-f\|_{\mathcal{L}^2(\SS^d)}\geq\big\|f_n-\mathcal{P}_{2^{q}}(f_n)\big\|_{\mathcal{L}^2(\SS^d)}-\|f-\mathcal{P}_{2^{q}}(f)\|_{\mathcal{L}^2(\SS^d)}\gtrsim n^{-\frac12}\underline{h}^{k+\frac12}\|f\|_{\mathcal{L}^2(\SS^d)}.
\end{equation}
    
\end{proof}

\begin{remark}
The saturation bound in Theorem~\ref{thm_main} requires antipodal quasi-uniformity. The separation from antipodal pairs is essential if one seeks a positive lower bound valid for every nonzero target: if $\theta_i^*=-\theta_j^*$, then
\[
(\theta_j^*\cdot\eta)^k=\sigma_k(\theta_j^*\cdot\eta)+(-1)^k\sigma_k(\theta_i^*\cdot\eta),\qquad\eta\in\SS^d,
\]
so this nonzero polynomial is represented exactly. In this case $\underline{h}=0$, and the configuration-dependent lower bound correctly becomes trivial.

For antipodally quasi-uniform families, Theorem~\ref{thm_main} is stronger than the standard form of saturation: every nonzero function in $\mathcal{H}^s(\SS^d)$ is excluded from approximation at a rate faster than $\mathcal{O}(n^{-\frac{d+2k+1}{2d}})$. For arbitrary quasi-uniform points and a general domain $\Omega$, we conjecture the standard lower bound statement: there exists a function that cannot be approximated at any rate $\mathcal{O}(n^{-\frac{d+2k+1}{2d}-\epsilon})$, $\epsilon>0$.

\end{remark}

\begin{remark}
The half-sphere interpretation is valid only modulo a fixed polynomial correction. Let $H\subset\SS^d$ contain one representative of each antipodal pair, and replace each $\theta_j^*$ by the representative $\theta_j^\sharp\in\{\theta_j^*,-\theta_j^*\}\cap H$. Since
\[
(-1)^k\sigma_k(-\theta\cdot\eta)+\sigma_k(\theta\cdot\eta)=(\theta\cdot\eta)^k,
\]
changing the orientation of a ridge changes it by a polynomial in $\PP_k(\SS^d)$. Consequently,
\begin{equation*}
L_n^k\big(\{\theta_j^*\}_{j=1}^n,\SS^d\big)+\PP_k(\SS^d)=L_n^k\big(\{\theta_j^\sharp\}_{j=1}^n,\SS^d\big)+\PP_k(\SS^d).
\end{equation*}
Thus full-sphere and half-sphere parameterizations can be compared after adjoining a fixed finite-dimensional polynomial space. This does not identify the unaugmented spaces, does not imply that an arbitrary quasi-uniform set has $\underline{h}\simeq n^{-1/d}$, and does not extend the sharp $n$-dependent lower bound beyond center families with the required antipodal separation.
\end{remark}

\section{Conclusion}\label{sec_concl}

In this paper, we have established a configuration-dependent saturation theorem for linearized shallow ReLU$^k$ approximation on the sphere. For arbitrary parameter sets, the lower bound is expressed as $n^{-1/2}\underline{h}^{k+1/2}$. For antipodally quasi-uniform families, this becomes the sharp rate $n^{-\frac{d+2k+1}{2d}}$, matching the known upper bound at the critical regularity $r=\frac{d+2k+1}{2}$. Thus the advantage of these linearized spaces over degree-$k$ finite elements is inherently bounded, while the dependence on more general center configurations is quantified explicitly by $\underline{h}$.

Our saturation theorem aligns neural network approximation with classical methods such as polynomial, spline, wavelet, and kernel approximations, where saturation phenomena are fundamental and well-documented. This underscores a universal structural limitation governing the performance of approximation schemes, extending even to nonlinear, adaptive methods such as neural networks. Practically, our results caution against overly optimistic views of shallow neural networks' capabilities, suggesting that their expressiveness—though superior—is ultimately limited by an intrinsic regularity threshold.

Looking forward, this saturation perspective naturally raises several research directions. Future studies might explore whether the same saturation order holds on a general domain $\Omega\subset\RR^d$. Moreover, \cite{siegel2022high} showed that nonlinear shallow ReLU$^k$ approximation can achieve $\mathcal{O}(n^{-(k+1)})$ for very smooth functions. Whether its saturation order is $k+1$, however, remains open.

\bibliographystyle{abbrv}
\bibliography{ref}

\end{document}